\documentclass[12pt, reqno]{amsart}
\usepackage{amssymb, amsthm, mathrsfs}

\newtheorem{thm}{Theorem}[section]

\newtheorem{lemma}[thm]{Lemma}
\newtheorem{theorem}[thm]{Theorem}

\numberwithin{equation}{section}

\theoremstyle{definition}
\newtheorem{rem}[thm]{Remark}

\newcommand{\al}{\alpha}

\newcommand{\e}{\varepsilon}
\newcommand{\la}{\lambda}
\renewcommand{\phi}{\varphi}

\renewcommand{\d}{\partial}
\newcommand{\R}{{\mathbb R}}
\newcommand{\br}[1]{\left\langle #1 \right\rangle}
\newcommand{\Case}[1]{\noindent \textbf{Case #1.}}

\newcommand{\nequiv}{\mathrel{\setbox0\hbox{$\equiv$}%
                     \rlap{\hbox{$\equiv$}}\hbox to \wd0{\hfil $/$\hfil}}}

\renewcommand{\qed}{\rule{3mm}{3mm}}
\renewenvironment{proof}
    {\vspace{1mm}\noindent\textbf{Proof.}}
    {\hspace*{\fill} $\qed$\vspace{1mm}}

\begin{document}
\title[Instability of steady states]
{Instability of steady states for nonlinear wave and heat equations}
\author{Paschalis Karageorgis}
\author{Walter A. Strauss}
\address{School of Mathematics, Trinity College, Dublin 2, Ireland}
\email{pete@maths.tcd.ie}
\address{Department of Mathematics and Lefschetz Center for Dynamical Systems,
Brown University, Providence, RI 02912}
\email{wstrauss@math.brown.edu}
\thanks{The second author was supported in part by NSF grant DMS-0405066}



\begin{abstract}
We consider time-independent solutions of hyperbolic equations such as $\d_{tt}u -\Delta u= f(x,u)$ where $f$ is convex
in $u$. We prove that linear instability with a positive eigenfunction implies nonlinear instability. In some cases the
instability occurs as a blow up in finite time. We prove the same result for parabolic equations such as $\d_t u -\Delta
u= f(x,u)$. Then we treat several examples under very sharp conditions, including equations with potential terms and
equations with supercritical nonlinearities.
\end{abstract}
\maketitle

\section{Introduction}
Given a linear second-order elliptic differential operator $L$ whose coefficients are smooth and bounded, consider the
parabolic equation
\begin{equation}\label{peq}
\d_t u + Lu = f(x,u), \quad\quad x\in \R^n
\end{equation}
along with its hyperbolic analogue
\begin{equation}\label{heq}
\d_t^2 u + a\d_t u + Lu =  f(x,u), \quad\quad x\in \R^n,
\end{equation}
where $f$ is a nonlinear term and $a\in \R$ is arbitrary (possibly zero).  A very important step in understanding the
behavior of general solutions lies in understanding the qualitative properties of special types of solutions.  In this
paper, we focus on time-independent solutions, also known as steady states, and we address their stability properties in
the context of both \eqref{peq} and \eqref{heq}.  Our main goal is to provide sufficient conditions under which
linearized instability can be used to draw conclusions about nonlinear instability.

Before we turn to our main results, however, let us first introduce some assumptions on the time-independent solution
$\phi$ and the nonlinear term $f$.  We are going to assume that
\begin{itemize}
\item[(A1)]
the equation $L\phi=f(x,\phi)$ has a $\mathcal{C}^2$ solution $\phi$;
\item[(A2)]
the adjoint linearized operator $L^\ast - f_u(x,\phi)$ has a negative eigenvalue $-\sigma^2$ and a corresponding
eigenfunction $\chi \in L^1(\R^n)\cap L^2(\R^n)$ that is non-negative;
\item[(A3)]
both $f(x,\phi)$ and $f_u(x,\phi)$ are bounded;
\item[(A4)]
the nonlinear term $f(x,s)$ is convex in $s$ and is $\mathcal{C}^1$.
\end{itemize}
Here, (A2) is mostly meant to ensure the presence of a negative eigenvalue; our assertions about the eigenfunction $\chi$
are already implied by this fact under pretty general conditions. Also, note that we do {\it not} require the steady
state $\phi$ to be bounded.

A non-technical description of our main result is that $\phi$ is a nonlinearly unstable solution of both \eqref{peq} and
\eqref{heq} whenever (A1)-(A4) hold; namely, solutions which start out close to $\phi$ need not remain close to $\phi$
for all times.  When it comes to the parabolic case \eqref{peq}, our result applies to a wide class of initial data,
including all initial data for which
\begin{equation}
u(x,0)>\phi(x).
\end{equation}
To establish the instability of $\phi$, we show that the norm
\begin{equation}\label{norm}
||u(x,t) - \phi(x)||_{L^\infty(\R^n)}
\end{equation}
must either grow exponentially at all times or else blow up in finite time.  Our main result for the hyperbolic case
\eqref{heq} is almost identical. It too applies to a wide class of initial data, including all initial data for which
\begin{equation}
u(x,0)>\phi(x), \quad\quad \d_t u(x,0) > 0
\end{equation}
and it shows that the energy norm\footnote{In each case, our norm agrees with the one dictated by the local existence
theory.}
\begin{equation}\label{norm2}
||u-\phi||_e \equiv ||u(x,t) - \phi(x)||_{H^1(\R^n)} + ||\d_t u(x,t)||_{L^2(\R^n)}
\end{equation}
must either grow exponentially at all times or else blow up in finite time.

Note that the above results prove instability in a sense that is much stronger, and perhaps more natural, than the usual
one. Namely, not only do they show that the solution exits any given neighborhood of $\phi$ in finite time, but they also
ensure that the solution does not reenter that neighborhood at any later time.

If one is willing to impose an additional positivity condition on the nonlinear term, then the above results can be
further improved to show that instability occurs by blow up.  This is the case, for instance, if one additionally assumes
that
\begin{itemize}
\item[(A5)]
there exist $C_0>0$ and $p>1$ such that $f(x,s) \geq C_0|s|^p$ for all $(x,s)\in\R^n\times \R$;
\item[(A6)]
the product $\phi\chi$ is integrable, where $\chi$ denotes the eigenfunction from (A2).
\end{itemize}
In our main results, we shall prove \textit{instability} assuming (A1)-(A4), and also \textit{instability by blow up}
assuming (A1)-(A6). As the reader should note, (A3)-(A6) may be safely ignored in the special case $f(u)= |u|^p$ for some
$p>1$, provided that $\phi$ is bounded.  For that special case, in particular, we establish instability by blow up
assuming (A1)-(A2) only.

To a large extent, our results are complementary to the abstract result of Grillakis, Shatah and Strauss \cite{GSS}.
These authors deal with arbitrary nonlinearities and the more general class of bound states, but their assumptions are
more restrictive than ours.  As we already mentioned, another advantage of our approach is that our notion of instability
is much stronger than the one used in \cite{GSS}. As we indicate later in the introduction, a final advantage is that our
result applies to a very sharp class of perturbations, which is not the case with the result of \cite{GSS}.

Section \ref{main} is devoted to the proof of our two main instability results, Theorems \ref{wave} and \ref{heat}. Our
proof is quite elementary and based on a variant of Kaplan's eigenfunction method \cite{Kap} which allows us to reduce
our analysis to the study of a certain functional. Although there are several other variants of this method, they are
concerned with the associated Cauchy problem rather than the instability problem; see \cite{Fu, Gl1, Ka2, LN, Lev, YZ1},
for instance. In the last three sections of this paper, we apply our two main results to treat some examples that we
briefly discuss below.

Our first and most general example appears in Theorem \ref{new}.  Here, we shall only describe a very special case of the
theorem that is itself of  independent interest.  Consider the nonlinear heat and wave equations with
potential\footnote{In the rest of the introduction, the possible damping term $a\d_t u$ is suppressed merely for the sake
of exposition.}
\begin{equation}\label{ex1}
\d_t^i u - \Delta u + V(x)\cdot u = |u|^p, \quad\quad x\in\R^n; \quad\quad i=1,2.
\end{equation}
We take $p>1$ and assume that $V$ is bounded, continuous and non-negative, but we make no other assumptions. Then Theorem
\ref{new} implies the instability of all non-negative, $H^1$ steady states that vanish at infinity.  A sufficient
condition on $V$ which ensures the existence of such steady states (for some $p$) is provided by \cite{SZ} and stated in
Theorem \ref{exam}.  Together with Theorem \ref{new}, these results impose no restrictions on the linearized operator; we
are not aware of any other general results with this feature.

In our second example, Theorem \ref{exp}, we focus on the two-dimensional equations\footnote{This is our only example for
which the space dimension is restricted.}
\begin{equation}\label{ex2}
\d_t^i u - \Delta u = e^u, \quad\quad x\in \R^2; \quad\quad i=1,2.
\end{equation}
The classification of all $\mathcal{C}^2$ steady states for which the nonlinear term is integrable is provided by
\cite{CL} and stated in Lemma \ref{expex}. They are {\it unbounded} and, as far as we know, their instability was not
previously known.  Theorem \ref{exp} shows they are all unstable.

In our last example, Theorem \ref{old}, we focus on the equations
\begin{equation}\label{ex3}
\d_t^i u - \Delta u = |u|^p, \quad\quad x\in \R^n; \quad\quad i=1,2,
\end{equation}
where $n>2$ and $p\geq\frac{n+2}{n-2}$.  When it comes to the existence of positive, $\mathcal{C}^2$, radially symmetric
steady states, the known results \cite{CGS, CL, DN, Li} are summarized in Theorem \ref{ex}.  Except for the special case
$p=\frac{n+2}{n-2}$, these steady states are not in $H^1$, so our first example is no longer applicable.  In section
\ref{LE}, we analyze the spectrum of the linearized operator and we find the exact values of $p$ for which a negative
eigenvalue emerges.  Using our main results, we then show that there is a critical value $p_c$, depending on the space
dimension, such that the steady states are \textit{nonlinearly unstable} if $p<p_c$ and \textit{linearly stable} if
$p\geq p_c$; see Theorem \ref{old}.

Our conclusions for equation \eqref{ex3} are new only in the hyperbolic case.  As for the parabolic equation, the
critical value $p_c$  emerged in \cite{GNW}, where a quite different approach was used. The instability result of
\cite{GNW} is weaker than ours because it applies to a smaller class of initial data, but the stability result given
there is stronger since it proves nonlinear stability instead. On the other hand, the approach in \cite{GNW} relies on
the maximum principle, so it cannot be applied to yield analogous results for the wave equation; see Remark \ref{re-h}
for more comments.

Finally, our result for the wave equation \eqref{ex3} is closely related to a recent paper of Krieger and Schlag
\cite{KS}.  These authors focus on the three-dimensional quintic case
\begin{equation}\label{ex4}
\d_t^2 u - \Delta u = |u|^5, \quad\quad x\in\R^3
\end{equation}
and a particular positive, $\mathcal{C}^2$ steady state $\phi$.  Theorem \ref{old} implies that $\phi$ is unstable, while
the result of \cite{KS} asserts the existence of a stable manifold associated with $\phi$.  Based on numerical evidence
from \cite{Num1, Num2}, this stable manifold ought to separate the \textit{small} perturbations for which solutions exist
globally from those for which solutions blow up.  In fact, Theorem \ref{old} provides a partial proof of this conjecture,
namely the blow up for all perturbations that lie strictly above the tangent plane to the stable manifold at the origin;
see Remark \ref{re-w} for more details.


\section{The main instability results}\label{main}
In this section, we prove our two main instability results regarding the solutions of a general elliptic equation of the
form
\begin{equation*}
L\phi = f(x,\phi), \quad\quad x\in\R^n.
\end{equation*}
Recall that $L$ is a linear, second-order elliptic differential operator whose coefficients are smooth and bounded, while
$f$ is a nonlinear term.  Our precise assumptions (A1)-(A6) were already mentioned in the introduction, so we shall not
bother to repeat them here.

First, we deal with the hyperbolic case and thus focus on the equation
\begin{equation}\label{1}
\d_t^2 u +a\d_t u + L u = f(x,u), \quad\quad u(x,0)= \phi(x)+\psi_0(x), \quad\quad \d_t u(x,0)= \psi_1(x).
\end{equation}
Since the steady state $\phi$ is an exact solution when $\psi_0\equiv \psi_1\equiv 0$, we are mostly concerned with the
case that the perturbation $(\psi_0,\psi_1)$ is small in some sense.  In our next result, we shall deal with
finite-energy perturbations.  Note, however, that $\phi$ itself need not be of finite-energy.

\begin{theorem}\label{wave} {\bf (Hyperbolic Equation)}
Let $a\in\R$.  Assume (A1)-(A4) and let $(\psi_0,\psi_1)\in H^1(\R^n)\times L^2(\R^n)$ be such that
\begin{equation}\label{per1}
\frac{a+\sqrt{a^2+4\sigma^2}}{2} \;\int_{\R^n} \chi(x)\,\psi_0(x) \:dx + \int_{\R^n} \chi(x)\,\psi_1(x) \:dx >0.
\end{equation}
Let $0<T\leq \infty$ and let $u$ be a solution of \eqref{1} on $[0,T)$ such that $u-\phi$ is continuous in $t$ with
values in the energy space and $f(x,u)$ is locally integrable.
\begin{itemize}
\item[(a)]
If $T=\infty$, then the energy norm
\begin{equation}\label{2}
||u(t)-\phi||_e \equiv ||u(\cdot,t)-\phi(\cdot)||_{H^1(\R^n)} + ||\d_t u(\cdot,t)||_{L^2(\R^n)}
\end{equation}
must grow exponentially.
\item[(b)]
Assume also (A5)-(A6).  Then $T<\infty$.
\end{itemize}
\end{theorem}

\begin{proof}
Consider the function
\begin{equation}\label{G}
G(t) = \int_{\R^n} \chi(x) \cdot w(x,t) \:dx, \quad\quad w(x,t)= u(x,t)- \phi(x).
\end{equation}
By (A2), one certainly has
\begin{equation*}
|G(t)| \leq ||\chi||_{L^2(\R^n)} \cdot ||w(t)||_{L^2(\R^n)} \leq C||w(t)||_e .
\end{equation*}
 Thus $G(t)$ is well-defined and bounded as long as the energy remains bounded; the energy
grows exponentially, provided that $G(t)$ does; and the energy becomes infinite whenever $G(t)$ does. In view of these
facts, our assertions about the energy \eqref{2} will follow once we have established analogous assertions for $G(t)$.

Let us first focus on part (a) and assume that  $w=u-\phi$ is continuous on $[0,\infty)$ with values in the energy space.
Then $w$ is a solution of
\begin{equation*}
\d_t^2 w +a\d_t w + Lw = f(x,w+\phi) -f(x,\phi)
\end{equation*}
in the sense of distributions. In view of our convexity assumption (A4), we have
\begin{equation*}
\d_t^2 w +a\d_t w + [L-f_u(x,\phi)] \,w \geq 0 .
\end{equation*}
Since $\chi(x)\ge0$, we may multiply the inequality by $\chi(x)\theta(t)$, where $\theta(t)$ is an arbitrary non-negative
test function, and integrate to obtain
\begin{align}\label{weak2}
-\int_0^t \int_{\R^n} \chi(x) \cdot \d_t w \cdot \theta'(\tau) &\:dx\,d\tau +
\int_0^t \int_{\R^n} \chi(x) \cdot a\d_t w \cdot \theta(\tau) \:dx\,d\tau \notag\\
&+ \int_0^t \int_{\R^n} [L^\ast - f_u(x,\phi)] \,\chi(x) \cdot w \cdot \theta(\tau) \:dx\,d\tau \geq 0 .
\end{align}
 To simplify the first two integrals, we note that
\begin{equation*}
G'(t) = \int_{\R^n} \chi(x)\cdot \d_t w(x,t) \:dx
\end{equation*}
is a continuous function of $t$ by our definition \eqref{G} because $\d_t w$ is continuous with values in $L^2(\R^n)$.
To simplify the third integral in \eqref{weak2}, we note that
\begin{equation*}
[L^\ast - f_u(x,\phi)] \,\chi = -\sigma^2 \chi
\end{equation*}
by our assumption (A2).  Thus, equation \eqref{weak2} reduces to
\begin{align*}
-\int_0^t G'(\tau) \cdot \theta'(\tau) &\:d\tau + a\int_0^t G'(\tau) \cdot \theta(\tau) \:d\tau -\sigma^2 \int_0^t
G(\tau) \cdot \theta(\tau) \:d\tau \geq 0
\end{align*}
for all non-negative test functions $\theta(t)$, or equivalently,
\begin{equation}\label{weak4}
G''(t) + aG'(t) -\sigma^2 G(t) \geq 0
\end{equation}
in the sense of distributions.  Next, we note that our assumption \eqref{per1} reads
\begin{equation*}
\frac{a+\sqrt{a^2+4\sigma^2}}{2} \cdot G(0) + G'(0) >0.
\end{equation*}
According to Lemma \ref{ODE1} below, both $G(t)$ and $G'(t)$ must then grow exponentially fast, so the proof of part (a)
is complete.

Next, we turn to part (b).  Suppose that $T=\infty$. As we have just shown, both $G(t)$ and $G'(t)$ must grow
exponentially fast. Using our assumptions (A3) and (A5), we have
\begin{align*}
\d_t^2 w +a\d_t w + [L-f_u(x,\phi)] \,w
&= f(x,w+\phi) -f(x,\phi) -f_u(x,\phi)w  \\
&\geq   C_0 |w+\phi|^p - C_1 - C_1|w|.
\end{align*}
Multiplying by the non-negative eigenfunction $\chi$ and integrating over space, we then easily find that
\begin{align}\label{eq22}
G''(t) +aG'(t) -\sigma^2 G(t) \geq C_0\int_{\R^n} \chi |w+\phi|^p \:dx - C_1\int_{\R^n} \chi \:dx -C_1\int \chi |w|\:dx,
\end{align}
the second derivative being understood in the sense of distributions. This can be justified, as above, by introducing a
non-negative test function $\theta(t)$ to avoid the second-order derivative.  Now, the last integral is at most
\begin{equation*}
\int \chi |w| \:dx \leq \int \chi|w+\phi|\:dx + \int \chi|\phi| \:dx
\end{equation*}
by the triangle inequality.  Since $\chi\in L^1$ by (A2) and $\chi\phi\in L^1$ by (A6), we deduce that
\begin{align*}
G''(t) +aG'(t) \geq \sigma^2 G(t) +C_0\int_{\R^n} \chi |w+\phi|^p \:dx - C_1\int \chi |w+\phi| \:dx - C_2 .
\end{align*}
Since $G(t)$ grows exponentially fast, this implies
\begin{align}\label{eq3}
G''(t) + aG'(t) &\geq C_0\int_{\R^n} \chi |w+\phi|^p \:dx - C_1 \int_{\R^n} \chi |w+\phi| \:dx \notag\\
&\equiv C_0A(t) - C_1B(t)
\end{align}
in the sense of distributions for all large enough $t$.

Now $B(t)$ itself grows exponentially fast, as the triangle inequality gives
\begin{equation}\label{eq4}
G(t) \leq \int_{\R^n} \chi |w|\:dx \leq \int_{\R^n} \chi |w+\phi|\:dx +\int_{\R^n} \chi |\phi|\:dx = B(t) +C_3.
\end{equation}
Similarly, $A(t)$ grows exponentially fast since H\"older's inequality gives
\begin{equation}
B(t)  \leq \left( \int_{\R^n} \chi \:dx \right)^{\frac{p-1}{p}} \left( \int_{\R^n} \chi |w+\phi|^p\:dx
\right)^{\frac{1}{p}} = C_4 A(t)^{1/p}.
\end{equation}
This actually forces $A(t)$ to grow faster than $B(t)$, namely
\begin{equation}\label{eq5}
A(t) \geq C_4^{-p} B(t)^p,
\end{equation}
hence $A(t)$ will eventually dominate $B(t)$.  Combining \eqref{eq3}, \eqref{eq4} and \eqref{eq5}, we now get
\begin{equation}\label{eq6}
G''(t) +aG'(t) \geq C_5 A(t) \geq C_6 B(t)^p \geq C_7 G(t)^p
\end{equation}
in the sense of distributions for all large enough $t$.  Moreover, both $G(t)$ and $G'(t)$ are eventually positive by
above.  Invoking Lemma \ref{ODE2} below, we reach the contradiction $T<\infty$.
\end{proof}

Next, we modify our previous approach to obtain a simple parabolic analogue of Theorem \ref{wave}.

\begin{theorem}\label{heat} {\bf (Parabolic Equation)}
Assume (A1)-(A4) and let $\psi_0\in L^\infty(\R^n)$ be continuous with
\begin{equation}\label{per2}
\int_{\R^n} \chi(x) \,\psi_0(x) \:dx > 0.
\end{equation}
Let $0<T\leq \infty$ and let $u$ be a solution of
\begin{equation}\label{3}
\d_t u+ Lu = f(x,u), \quad\quad u(x,0)= \phi(x)+\psi_0(x)
\end{equation}
on $[0,T)$ such that $u-\phi$ is continuous and bounded.
\begin{itemize}
\item[(a)]
If $T=\infty$, then the norm $||u-\phi||_{L^\infty(\R^n)}$ must grow exponentially.
\item[(b)]
Assume also (A5)-(A6). Then $T<\infty$.
\end{itemize}
\end{theorem}

\begin{proof}
Once again, it suffices to establish analogous assertions for the function
\begin{equation*}
G(t) = \int_{\R^n} \chi(x) \cdot w(x,t) \:dx, \quad\quad w(x,t)= u(x,t)- \phi(x).
\end{equation*}
Since
\begin{equation*}
|G(t)| \leq ||\chi||_{L^1(\R^n)} \cdot ||w||_{L^\infty(\R^n)},
\end{equation*}
the $L^\infty$-norm has to either grow exponentially or blow up whenever $G(t)$ does.

We now apply the same argument with minor changes.  Note that $w=u-\phi$ satisfies
\begin{equation*}
\d_t w + Lw = f(x,w+\phi)- f(x,\phi), \quad\quad w(x,0)= \psi_0(x)
\end{equation*}
in the sense of distributions.  Arguing as before, we then get the estimate
\begin{equation}\label{eqq}
G'(t) -\sigma^2 G(t) \geq 0
\end{equation}
instead of \eqref{weak4}.  Since $G(0)>0$ by assumption, the exponential growth of $G(t)$ follows directly.

Under the additional assumptions (A5)-(A6), our previous approach yields the estimate
\begin{align*}
G'(t) -\sigma^2 G(t) \geq C_0\int_{\R^n} \chi |w+\phi|^p \:dx - C_1\int_{\R^n} \chi \:dx -C_1\int \chi |w|\:dx
\end{align*}
instead of \eqref{eq22}.  Then the remaining part of our proof applies verbatim to give
\begin{equation*}
G'(t) \geq C_7 G(t)^p
\end{equation*}
instead of \eqref{eq6}, in the sense of distributions for all large enough $t$.   Since $G(t)$ is positive, it is easy to
deduce that $T<\infty$, as needed.
\end{proof}

\begin{lemma}\label{ODE1}
Let $a\in \R$ and $b>0$.  Suppose $y(t)$ is a $\mathcal{C}^1$ function such that
\begin{equation*}
y'' + ay' - by \geq 0
\end{equation*}
on some interval $[0,T)$ in the sense of distributions.  If
\begin{equation}\label{ic}
\frac{a+\sqrt{a^2+4b}}{2} \cdot y(0) + y'(0) > 0,
\end{equation}
then both $y(t)$ and $y'(t)$ must grow exponentially on $[0,T)$.
\end{lemma}

\begin{proof}
Let $\la_1<0<\la_2$ be the roots of the characteristic equation $\la^2+ a\la- b=0$ and set $z= y'-\la_1 y$.  Then
\begin{align*}
z'-\la_2z =y''+ay'-by\geq 0.
\end{align*}
Using the test function $\exp(-\lambda_2t)\,\theta(t)$, where $\theta(t)$ is another test function, it follows that
\begin{equation*}
-\int_0^t z(\tau) \exp (-\lambda_2\tau) \,\theta'(\tau) \:d\tau \geq 0.
\end{equation*}
Choosing $\theta(t)$ to be an approximation of the characteristic function of the interval $(0,t)$, we easily deduce that
$z(t) \geq e^{\lambda_2t} z(0)$. Thus $y'-\la_1 y\geq e^{\la_2 t} z(0)$, which implies that
\begin{equation*}
y(t) \geq e^{\la_1t} y(0) + \frac{e^{\la_2 t}-e^{\la_1t}}{\la_2-\la_1} \cdot z(0).
\end{equation*}
Since $\la_1 < 0< \la_2$ by above, the exponential growth of $y$ then follows, provided that
\begin{equation*}
z(0) = y'(0) -\la_1 y(0) = y'(0) + \frac{a+\sqrt{a^2+4b}}{2} \cdot y(0)
\end{equation*}
is positive.  In view of our assumption \eqref{ic}, the exponential growth of $y$ thus follows.

Next, we use a similar argument to find that $w= y'-\la_2y$ satisfies
\begin{align*}
w'-\la_1w =y''+ay'-by\geq 0,
\end{align*}
hence also $w(t)\geq e^{\la_1 t} w(0)$.  Since $y(t)$ grows exponentially by above, the equation
\begin{equation*}
y'(t) \geq \la_2y(t) + e^{\la_1t} w(0)
\end{equation*}
then forces $y'(t)$ to grow exponentially fast as well because $\la_1<0<\la_2$.
\end{proof}

\begin{lemma}\label{ODE2}
Let $a\in \R$, $b>0$ and $p>1$.  Suppose $y(t)$ is a non-negative $\mathcal{C}^1$ function such that
\begin{equation*}
y(T_1)> 0, \quad\quad y'(T_1)> 0, \quad\quad y''(t) + ay'(t) \geq by(t)^p
\end{equation*}
on some interval $[T_1,T_2)$ in the sense of distributions.  Then $T_2<\infty$.
\end{lemma}

\begin{proof}
Suppose first that $y\in \mathcal{C}^2$.  Then the case $a=1$ is treated for instance in Proposition 3.1 of \cite{TY};
the case $a>0$ is quite similar; and the case $a\leq 0$ is much easier.  If $y$ is merely $\mathcal{C}^1$, one can simply
repeat the same argument using test functions as in the preceding lemma; we omit the details.
\end{proof}

\section{Convex nonlinearity with potential term}
In this section, we apply our main results to study non-negative $H^1$ solutions of the equation
\begin{equation}\label{neweq}
-\Delta \phi +V(x)\cdot \phi= f(\phi), \quad\quad x\in\R^n.
\end{equation}
Although our approach applies verbatim to the more general case $f(x,\phi)$, the assumptions of our next result would
have to be suitably modified in that case; we shall only deal with \eqref{neweq} for the sake of simplicity.

\begin{theorem}\label{new}
Assume that the following conditions hold:
\begin{itemize}
\item[(B1)]
$V$ is continuous and bounded on $\R^n$;
\item[(B2)]
the essential spectrum of $-\Delta+V$ is contained in $[0,\infty)$;
\item[(B3)]
equation \eqref{neweq} has a non-negative $\mathcal{C}^2$ solution $\phi\in H^1$ that vanishes at infinity;
\item[(B4)]
$f\in \mathcal{C}^1(\R)$ is convex with $f(0)=f'(0)=0$ and $f(\phi)$ is not identically zero.
\end{itemize}
Then $\phi$ is an unstable solution of any of the equations
\begin{equation*}
\d_t u + (-\Delta+V)u = f(u), \quad\quad \d_t^2 u + a\d_t u + (-\Delta+V)u = f(u),
\end{equation*}
where $a\in\R$.  More precisely, the conclusions of Theorems \ref{wave}(a) and \ref{heat}(a) remain valid.
\end{theorem}

\begin{rem}
By Weyl's theorem, our spectral assumption (B2) essentially requires that $V(x)$ not be negative as $|x|\to\infty$.   We
do not know of any existence results regarding \eqref{neweq} for which this assumption is violated; there are several
existence results \cite{BL, SZ, ST-ex} in case it holds.
\end{rem}

\begin{proof}
We verify the assumptions (A1)-(A4) of Theorems \ref{wave} and \ref{heat}.  Note that the existence assumption (A1) holds
by (B3), while the convexity assumption (A4) holds by (B4).  Since $\phi$ is bounded by (B3), both $f(\phi)$ and
$f'(\phi)$ are bounded as well, so the assumption (A3) also holds. To check the remaining assumption (A2), we consider
the self-adjoint linearized operator
\begin{equation*}
\mathscr{L} = -\Delta +V -f'(\phi).
\end{equation*}
Noting that $\phi$ is a solution of \eqref{neweq}, we use (B4) to find that
\begin{equation*}
\mathscr{L\phi} = f(\phi) - \phi \cdot f'(\phi) = f(\phi) - f(0) - \phi \cdot f'(\phi) \leq 0.
\end{equation*}
Let $\br{\,,}$ denote the standard inner product on $L^2(\R^n)$.  Then the expression
\begin{equation*}
\br{\mathscr{L}\phi,\phi} = \int_{\R^n} |\nabla \phi|^2\:dx + \int_{\R^n} [V  -  f'(\phi)]\cdot \phi^2\:dx
\end{equation*}
is non-positive because $\mathscr{L}\phi\leq 0 \leq \phi$ and is finite because $\phi\in H^1$ and $V-f'(\phi)\in
L^\infty$. Thus, the first eigenvalue of $\mathscr{L}$ must also be non-positive.  It can only be zero if
\begin{equation*}
\mathscr{L} \phi = f(\phi) - \phi\cdot f'(\phi)
\end{equation*}
is identically zero, in which case $f(y)= yf'(y)$ for all $y$ in the image of $\phi$.  Solving this ordinary differential
equation, we get $f(y)=Cy$ for some constant $C$.   Then the fact that $f'(0)=0$ implies $f(y)= 0$ for all $y$ in the
image of $\phi$.  Since this violates our assumption (B4), however, the first eigenvalue must actually be negative.

Now that the first eigenvalue is known to be negative, the remaining assertions of (A2) follow by standard facts.  Using
a variational argument, it is well-known that the first eigenfunction can be chosen to be positive (see section 11 in
\cite{LL}, for instance).  Finally, the first eigenfunction is in $L^1\cap L^2$ because it decays exponentially by
Agmon's estimate; see Theorem C.3.5 in \cite{Sim1} or else a more general result of Nakamura \cite{Nak}. In fact, Agmon's
estimate applies to any eigenfunction for which the associated eigenvalue lies below the bottom of the essential
spectrum.  Since
\begin{equation*}
\lim_{|x|\to \infty} f'(\phi(x)) = f'(0) = 0
\end{equation*}
by (B3) and (B4), an application of Weyl's theorem (see Theorem 1.1 in \cite{Sim3}) gives
\begin{equation*}
\sigma_{\text{ess}}(\mathscr{L})= \sigma_{\text{ess}}(-\Delta+V) \subset [0,\infty)
\end{equation*}
because of (B2).  Thus, the first eigenfunction must decay exponentially, as needed.
\end{proof}

In our next result, we give a typical application of Theorem \ref{new}.  Our precise assumptions on the potential $V(x)$
are taken from \cite{SZ}; they are merely meant to ensure that equation \eqref{szeq} has a non-negative $\mathcal{C}^2$
solution $\phi\in H^1$ that vanishes at infinity.  Needless to say, Theorem \ref{new} implies the instability of such
solutions for a much wider class of potentials.  The only advantage of our next result is that all the assumptions are
verified explicitly.

\begin{theorem}\label{exam}
Let $n>2$ and $1<p<\frac{n+2}{n-2}$.  Assume that $V(x)$ is radially symmetric, continuous and locally H\"older
continuous with
\begin{equation}\label{szv}
C_1(1+|x|)^{-l} \leq V(x) \leq C_2
\end{equation}
for some constants $C_1,C_2>0$ and some $0\leq l<\frac{2(n-1)(p-1)}{p+3}$. Then the equation
\begin{equation}\label{szeq}
-\Delta \phi +V(x)\cdot \phi= \phi^p, \quad\quad x\in\R^n
\end{equation}
has a positive, $\mathcal{C}^2$, radially symmetric and exponentially decaying solution $\phi$.  Moreover, $\phi$ is an
unstable solution of any of the equations
\begin{equation*}
\d_t u + (-\Delta+V)u = |u|^p, \quad\quad \d_t^2 u + a\d_t u + (-\Delta+V)u = |u|^p,
\end{equation*}
where $a\in\R$.  More precisely, the conclusions of Theorems \ref{wave} and \ref{heat} (both (a) and (b)) remain valid.
That is, the instability occurs by blow up in finite time.
\end{theorem}

\begin{proof}
Under the given hypotheses, our assertions about the elliptic equation \eqref{szeq} follow from Theorem 1.2 in \cite{SZ}.
These ensure that the assumption (B3) from the previous theorem holds with $f(u)= |u|^p$.  It is clear that this
nonlinear term satisfies (B4). Moreover, $V(x)$ is bounded and non-negative by \eqref{szv}, so the remaining assumptions
(B1)-(B2) hold as well. According to the previous theorem then, the conclusions of Theorems \ref{wave} and \ref{heat}
remain valid, indeed.

To show that instability occurs by blow up, it remains to check that the additional assumptions (A5)-(A6) of these
theorems remain valid as well.  Since $f(u)=|u|^p$ satisfies (A5) and $\phi$ is bounded, these two assumptions hold
trivially and the proof is complete.
\end{proof}


\section{Exponential nonlinearity in 2D}
In this section, we apply our main results to study solutions of the equation
\begin{equation}\label{expeq}
- \Delta \phi = e^\phi, \quad\quad e^\phi\in L^1(\R^2).
\end{equation}
The classification of all $\mathcal{C}^2$ solutions is due to Chen and Li \cite{CL}. It is perhaps worth noting that
these solutions are {\it unbounded} and not necessarily of one sign.

\begin{lemma}\label{expex}
Every $\mathcal{C}^2$ solution of \eqref{expeq} is of the form
\begin{equation*}
\phi(x) = \log \Bigl[ 32\la^2 \cdot (4+\la^2|x-y|^2)^{-2} \Bigr]
\end{equation*}
for some $\la>0$ and some $y\in\R^2$.
\end{lemma}

\begin{theorem}\label{exp}
Every $\mathcal{C}^2$ solution of \eqref{expeq} is an unstable solution of any of the equations
\begin{equation*}
\d_t u -\Delta u = e^u, \quad\quad \d_t^2 u + a\d_t u -\Delta u = e^u,
\end{equation*}
where $a\in\R$.  More precisely, the conclusions of Theorems \ref{wave}(a) and \ref{heat}(a) remain valid.
\end{theorem}

\begin{proof}
We set $f(u)=e^u$ and verify the assumptions (A1)-(A4) of Theorems \ref{wave} and \ref{heat}.  The existence assumption
(A1) and the convexity assumption (A4) obviously hold in this case.  Assumption  (A3) holds  as well because
\begin{equation*}
f(\phi) = f'(\phi) = 32\la^2 \cdot (4+\la^2|x-y|^2)^{-2}
\end{equation*}
is bounded.  To show that (A2) holds as well, we focus on the linearized operator
\begin{equation*}
\mathscr{L} = - \Delta - \exp\phi(x)= -\Delta -32\la^2 \cdot \left( 4+\la^2 |x-y|^2 \right)^{-2}.
\end{equation*}
Its essential spectrum is
\begin{equation*}
\sigma_{\text{ess}}(\mathscr{L}) = \sigma_{\text{ess}}(-\Delta) = [0,\infty)
\end{equation*}
by Weyl's theorem since $f'(\phi)$ is bounded and vanishes at infinity. Exactly as in the proof of Theorem \ref{new}, the
first eigenfunction will  be positive and exponentially decaying, so long as the first eigenvalue is negative.  In
particular, it suffices to show that the associated energy
\begin{equation}\label{en}
E(\zeta) = \int_{\R^2} |\nabla \zeta (x)|^2 \:dx - \int_{\R^2} \exp\phi(x)\cdot \zeta(x)^2 \:dx
\end{equation}
is negative for some test function $\zeta\in H^1(\R^2)$.   Choosing
\begin{equation*}
\zeta(x) = \left( 4+\la^2 |x-y|^2 \right)^{-2},
\end{equation*}
the  energy \eqref{en} is given by
\begin{align*}
E(\zeta) &= 16\la^4 \int_{\R^2} \frac{|x-y|^2\:\:dx}{(4+\la^2 |x-y|^2)^6} - 32\la^2 \int_{\R^2}
\frac{dx}{(4+\la^2 |x-y|^2)^6} \\
&= 32\la^4\pi \int_0^\infty \frac{r^3\:dr}{(4+\la^2 r^2)^6} - 64\la^2\pi \int_0^\infty
\frac{r\,dr}{(4+\la^2 r^2)^6}
\end{align*}
and then a short computation gives
\begin{align*}
E(\zeta) &= 16\pi \int_4^\infty (s-4) s^{-6} \:ds - 32\pi \int_4^\infty s^{-6} \:ds = -\frac{\pi}{320} \,.
\end{align*}
This implies the existence of a negative eigenvalue and completes the proof.
\end{proof}


\section{Power nonlinearity with zero mass}\label{LE}
In this last section, we focus on positive solutions of the equation
\begin{equation}\label{ee1}
-\Delta \phi (x)= \phi(x)^p, \quad\quad x\in \R^n.
\end{equation}
Before we state our results,  it is convenient to introduce the quadratic polynomial
\begin{equation}\label{Q1}
Q(\al) = \al(n-2-\al).
\end{equation}
This polynomial arises naturally through the computation
\begin{equation*}
-\Delta |x|^{-\al} = \al(n-2-\al)\cdot |x|^{-\al-2}
\end{equation*}
and it is closely related to Hardy's inequality
\begin{equation}\label{hl}
\int_{\R^n} |\nabla u(x)|^2 \:dx \geq \left( \frac{n-2}{2} \right)^2 \int_{\R^n} |x|^{-2} u(x)^2\:dx,
\end{equation}
which is valid for all $u\in H^1(\R^n)$ and each $n > 2$.  Namely, the coefficient on the right side of \eqref{hl} is the
maximum value of $Q(\al)$ and it is known to be sharp in the following sense.

\begin{lemma}\label{ene}
Let $n>2$ and let $V$ be a bounded function on $\R^n$ which vanishes at infinity.  If there exists some $\e>0$ such that
\begin{equation*}
V(x) \leq -(1+\e)\cdot \left( \frac{n-2}{2} \right)^2 \cdot |x|^{-2}
\end{equation*}
for all large enough $|x|$, then the operator $-\Delta+V$ has infinite negative spectrum.
\end{lemma}

For a proof of Hardy's inequality \eqref{hl}, see pg.~169 in \cite{RS2}.  For a proof of the  lemma, see the appendix in
\cite{DHS}, for instance.

We are now ready to state the known existence results for positive steady states. Parts (a) and (b) can be found in
either \cite{CGS} or \cite{CL}.  For part (c), see Theorem 1 in \cite{Li} and Theorem 5.26 in \cite{DN}.

\begin{theorem}\label{ex}
Let $n\geq 1$ and $p>1$.  Denote by $Q$ the quadratic in \eqref{Q1}.

\begin{itemize}
\item[(a)]
If either $n=1,2$ or $p< \frac{n+2}{n-2}$, then equation \eqref{ee1} has no positive, $\mathcal{C}^2$ solutions.

\item[(b)]
If $n>2$ and $p= \frac{n+2}{n-2}$, then any positive $\mathcal{C}^2$ solution of \eqref{ee1} is of the form
\begin{equation}\label{phila}
\phi_\la(x) = \left( \frac{\la \sqrt{n(n-2)}}{\la^2 + |x-y|^2} \right)^{\frac{2}{p-1}}
\end{equation}
for some $\la>0$ and some $y\in\R^n$.

\item[(c)]
If $n>2$ and $p>\frac{n+2}{n-2}$, then the positive, $\mathcal{C}^2$, radially symmetric solutions of \eqref{ee1} form a
one-parameter family $\{ \phi_\al \}_{\al>0}$, where each $\phi_\al$ satisfies
\begin{equation}\label{phia}
\phi_\al(0)= \al, \quad\quad \lim_{|x|\to \infty} |x|^2 \,\phi_\al(x)^{p-1} = Q\left( \frac{2}{p-1} \right) >0.
\end{equation}
\end{itemize}
\end{theorem}

In what follows, we use the previous two results to analyze the spectrum of the linearized operator associated with
\eqref{ee1}. In particular, we find the exact values of $p$ for which a negative eigenvalue emerges. Although a
preliminary characterization of these values is provided by the next two lemmas, a more concrete characterization will be
given in Theorem \ref{old}.

\begin{lemma}\label{ip}
Let $n>2$ and $p\geq \frac{n+2}{n-2}$.  Let $Q$ be the quadratic in \eqref{Q1} and let $\phi$ denote any one of the
steady states provided by Theorem \ref{ex}.  Then $-\Delta-p\phi^{p-1}$ has a negative eigenvalue if
\begin{equation}\label{ic1}
\left( \frac{n-2}{2} \right)^2 < p\cdot Q \left( \frac{2}{p-1} \right).
\end{equation}
In particular, it has a negative eigenvalue if $p=\frac{n+2}{n-2}$.
\end{lemma}

\begin{proof}
Suppose first that $p>\frac{n+2}{n-2}$.  Then \eqref{phia} and \eqref{ic1} allow us to find some $\e>0$ such that
\begin{equation*}
\lim_{|x|\to \infty} |x|^2 \,\phi(x)^{p-1} = Q\left( \frac{2}{p-1} \right) > p^{-1}(1+2\e) \cdot
\left( \frac{n-2}{2} \right)^2.
\end{equation*}
Thus
\begin{equation*}
V(x) = -p\phi(x)^{p-1} < -(1+\e)\cdot \left( \frac{n-2}{2} \right)^2 \cdot |x|^{-2}
\end{equation*}
for all large enough $|x|$, so the existence of a negative eigenvalue follows by Lemma \ref{ene}.

Suppose now that $p=\frac{n+2}{n-2}$.  Then inequality \eqref{ic1} automatically holds because
\begin{equation*}
\left( \frac{n-2}{2} \right)^2 = Q\left( \frac{n-2}{2} \right) = Q\left( \frac{2}{p-1} \right) < p\cdot Q\left(
\frac{2}{p-1} \right)
\end{equation*}
for this particular case.  According to part (b) of Theorem \ref{ex}, we also have
\begin{equation*}
- \Delta - p\phi(x)^{p-1}= -\Delta -\la^2 \,n(n+2) \cdot \left( \la^2 + |x-y|^2 \right)^{-2}
\end{equation*}
for some $\la>0$ and some $y\in \R^n$.  Thus it suffices to check that the associated energy
\begin{equation}\label{en2}
E(\zeta) = \int_{\R^n} |\nabla \zeta (x)|^2 \:dx - \int_{\R^n} p\phi(x)^{p-1} \:\zeta(x)^2 \:dx
\end{equation}
is negative for some test function $\zeta\in H^1(\R^n)$.  Let us now consider the test function
\begin{equation*}
\zeta(x) = \left( \la^2 + |x-y|^2 \right)^{-n/2-1},
\end{equation*}
whose energy \eqref{en2} is given by
\begin{equation*}
E(\zeta) = (n+2)^2 \int_{\R^n} \frac{|x-y|^2\:\:dx}{(\la^2+|x-y|^2)^{n+4}} - \la^2 \,n(n+2) \int_{\R^n}
\frac{dx}{(\la^2+|x-y|^2)^{n+4}} \,.
\end{equation*}
Except for a positive factor, this expression is equal to
\begin{equation}\label{en3}
\widetilde{E}(\zeta) = (n+2) \int_0^\infty \frac{r^{n+1}\:\:dr}{(\la^2+r^2)^{n+4}} - \la^2 \,n \int_0^\infty
\frac{r^{n-1}\:dr}{(\la^2+r^2)^{n+4}} \,.
\end{equation}
Moreover, an integration by parts gives
\begin{align*}
\int_0^\infty \frac{r^{n+1}\:\:dr}{(\la^2+r^2)^{n+4}}
&= \frac{n}{2(n+3)} \int_0^\infty \frac{r^{n-1}\:\:dr}{(\la^2+r^2)^{n+3}} \\
&= \frac{n}{2(n+3)} \left[ \la^2 \int_0^\infty \frac{r^{n-1}\:\:dr}{(\la^2+r^2)^{n+4}} + \int_0^\infty
\frac{r^{n+1}\:\:dr}{(\la^2+r^2)^{n+4}} \right] ;
\end{align*}
hence also
\begin{align*}
\int_0^\infty \frac{r^{n+1}\:\:dr}{(\la^2+r^2)^{n+4}} = \frac{\la^2 n}{n+6} \int_0^\infty
\frac{r^{n-1}\:\:dr}{(\la^2+r^2)^{n+4}} \,.
\end{align*}
Inserting this equality in \eqref{en3}, we arrive at
\begin{equation*}
\widetilde{E}(\zeta) = \la^2 n\left( \frac{n+2}{n+6} - 1\right) \int_0^\infty \frac{r^{n-1}\:dr}{(\la^2+r^2)^{n+4}} = -
\frac{4\la^2 n}{n+6} \int_0^\infty \frac{r^{n-1}\:dr}{(\la^2+r^2)^{n+4}} \,.
\end{equation*}
Since this expression is negative, the energy \eqref{en2} must also be negative, as needed.
\end{proof}

\begin{lemma}\label{sp}
Let $n>2$ and $p> \frac{n+2}{n-2}$.  Let $Q$ denote the quadratic in \eqref{Q1} and assume that
\begin{equation}\label{sc1}
\left( \frac{n-2}{2} \right)^2 \geq p\cdot Q \left( \frac{2}{p-1} \right).
\end{equation}
Then the steady states provided by Theorem \ref{ex} satisfy
\begin{equation}\label{bel}
|x|^2 \,\phi(x)^{p-1} \leq Q\left( \frac{2}{p-1} \right)
\end{equation}
for each $x\in\R^n$, and the operator $-\Delta-p\phi^{p-1}$ has no negative spectrum.
\end{lemma}

\begin{proof}
First suppose the estimate \eqref{bel} does hold.  Using our assumption \eqref{sc1}, we then get
\begin{equation*}
-p\phi(x)^{p-1} \geq -p \cdot Q\left( \frac{2}{p-1} \right) \cdot |x|^{-2} \geq -\left( \frac{n-2}{2} \right)^2 \cdot
|x|^{-2}
\end{equation*}
for each $x\in\R^n$.  According to Hardy's inequality \eqref{hl}, this already implies that
\begin{equation*}
\int_{\R^n} |\nabla u(x)|^2 \:dx - \int_{\R^n} p\phi(x)^{p-1}\cdot u(x)^2\:dx \geq 0
\end{equation*}
for all $u\in H^1(\R^n)$, so that the operator $-\Delta-p\phi^{p-1}$ has no negative spectrum, indeed.

Let us now prove \eqref{bel}.  Set $k=\frac{2}{p-1}$ and consider the function
\begin{equation*}
W(s) = e^{ks} \phi(e^s)= r^k \phi(r), \quad\quad s=\log r= \log |x|.
\end{equation*}
Then $W(s)$ is positive and bounded because of \eqref{phia}.  It also satisfies the equation
\begin{equation}\label{11}
Q(k-\d_s) W(s) = W(s)^p
\end{equation}
by Lemma \ref{cco} below. We note that $s$ ranges over $(-\infty,\infty)$ as $r$ ranges from $0$ to $\infty$, while
\begin{equation*}
\lim_{s\to -\infty} W(s) = \lim_{r\to 0^+} r^k\phi(r) = 0.
\end{equation*}
The derivatives of $W(s)$ must also vanish at $s=-\infty$ because
\begin{equation*}
\lim_{s\to -\infty} W'(s) = \lim_{r\to 0^+} r\cdot \d_r [r^k\phi(r)] = 0,
\end{equation*}
and so on.  Now the fact that $f(x)=x^p$ is convex on $(0,\infty)$ ensures that
\begin{equation*}
f(W(s)) - f\left( Q(k)^{\frac{1}{p-1}} \right) \geq f' \left( Q(k)^{\frac{1}{p-1}} \right) \cdot \left( W(s) -
Q(k)^{\frac{1}{p-1}} \right),
\end{equation*}
that is,
\begin{align*}
W(s)^p - Q(k)^{\frac{p}{p-1}} &\geq pQ(k)\cdot \left( W(s)-Q(k)^{\frac{1}{p-1}} \right).
\end{align*}
Inserting this inequality into \eqref{11}, we thus find
\begin{equation}\label{12}
Q(k-\d_s) W(s) - Q(k)^{\frac{p}{p-1}} \geq pQ(k)\cdot \left( W(s)-Q(k)^{\frac{1}{p-1}} \right).
\end{equation}
To eliminate the constant term on the left hand side, we change variables by
\begin{equation*}
Y(s) = W(s) - Q(k)^{\frac{1}{p-1}}.
\end{equation*}
Then equation \eqref{12} can also be written in the equivalent form
\begin{equation}\label{13}
\Bigl[ pQ(k) - Q(k-\d_s) \Bigr] \,Y(s) \leq 0.
\end{equation}

We will prove below that the characteristic polynomial
\begin{equation}\label{ch1}
\mathscr{P}(\la) \equiv pQ(k) -Q(k-\la)
\end{equation}
has two negative roots $\la_1,\la_2$. Temporarily assuming this and noting that  $\mathscr{P}(\la)$ has its highest-order
coefficient equal to $1$ by \eqref{Q1}, we can then factor the ordinary differential inequality \eqref{13} as
\begin{equation*}
Z(s) \equiv Y'(s)-\la_1 Y(s),\quad\quad Z'(s)-\la_2 Z(s)\leq 0.
\end{equation*}
The last inequality makes $e^{-\la_2s}Z(s)$ decreasing, so we must actually have
\begin{equation*}
e^{-\la_2s} Z(s) \leq \lim_{s\to -\infty} e^{-\la_2s} Z(s) = 0
\end{equation*}
because $\la_2<0$ by above.  This gives $Y'(s)-\la_1 Y(s)= Z(s)\leq 0$, and we similarly find
\begin{equation*}
e^{-\la_1s} Y(s) \leq \lim_{s\to -\infty} e^{-\la_1s} Y(s) = 0
\end{equation*}
because $\la_1<0$ as well.  Now that $Y(s)\leq 0$, we have
\begin{equation}\label{end}
W(s) \leq Q(k)^{\frac{1}{p-1}}.
\end{equation}
This is precisely the desired inequality \eqref{bel} since $W(s) = r^{\frac{2}{p-1}} \phi(r)$ by definition.

Thus it remains to show that the quadratic in \eqref{ch1} has two negative roots in case of assumption \eqref{sc1}. Since
the product of the two roots is
\begin{equation}\label{P0}
\mathscr{P}(0) = (p-1) \cdot Q(k) > 0
\end{equation}
by \eqref{ch1} and \eqref{phia}, it suffices to show that $\mathscr{P}(\la_*)\leq 0$ for some $\la_*<0$.  We choose
\begin{equation*}
\la_* = k- \frac{n-2}{2} = \frac{2}{p-1} - \frac{n-2}{2}
\end{equation*}
and note that $\la_*<0$ because $p>\frac{n+2}{n-2}$ by assumption. Combining \eqref{ch1}, \eqref{Q1} and \eqref{sc1}, we
find
\begin{equation*}
\mathscr{P}(\la_*) = pQ(k) - Q\left( \frac{n-2}{2} \right) = pQ(k) - \left( \frac{n-2}{2} \right)^2 \leq 0.
\end{equation*}
Thus $\mathscr{P}(\la)$ has two negative roots.
\end{proof}

\begin{lemma}\label{cco}
Let $n>2$.  Let $p\geq \frac{n+2}{n-2}$ and $k=\frac{2}{p-1}$. Given any one of the steady states $\phi$ provided by
Theorem \ref{ex}, the function
\begin{equation*}
W(s) = e^{ks} \phi(e^s) = r^k \phi(r), \quad\quad s=\log r = \log |x|
\end{equation*}
must then satisfy the ordinary differential equation
\begin{equation}\label{ne1}
(\d_s - k)(\d_s-k +n-2) \,W(s) = -W(s)^p.
\end{equation}
Moreover, this ordinary differential equation can be written in the form
\begin{equation*}
Q(k-\d_s)W(s) = W(s)^p,
\end{equation*}
where $Q$ denotes the quadratic in \eqref{Q1}.
\end{lemma}

\begin{proof}
Since $\d_r = e^{-s}\d_s$, a short computation allows us to write the radial Laplacian as
\begin{equation*}
\Delta = \d_r^2 + (n-1) r^{-1}\d_r = e^{-2s} (n-2+ \d_s)\d_s.
\end{equation*}
Since $-\Delta \phi = \phi^p$ and $kp= k+2$ by assumption, this easily leads us to
\begin{equation*}
-W^p = -e^{kps} \phi^p = e^{ks+2s} \Delta \phi = e^{ks} (n-2+ \d_s)\d_s (e^{-ks} W).
\end{equation*}
Using the operator identity $\d_s e^{-ks} = e^{-ks} (\d_s-k)$ twice, we then get
\begin{equation*}
-W^p = e^{ks} (n-2+ \d_s)e^{-ks}(\d_s-k) W = (n-2+ \d_s-k)(\d_s-k) W.
\end{equation*}
This proves our first assertion \eqref{ne1}, from which our second assertion follows trivially.
\end{proof}

\begin{theorem}\label{old}
Let $n>2$ and $p\geq \frac{n+2}{n-2}$.  Let $\phi$ denote any one of the steady states provided by Theorem \ref{ex} and
set
\begin{equation}\label{pc}
p_c = \left\{
\begin{array}{cc}
\infty &\text{if\: $n\leq 10$}\\
\frac{n^2-8n+4+8\sqrt{n-1}}{(n-2)(n-10)} &\text{if\: $n> 10$}\\
\end{array}\right\}
\end{equation}
for convenience.  Then $p_c>\frac{n+2}{n-2}$ and the following dichotomy holds.

If $p<p_c$, then $\phi$ is a nonlinearly unstable solution of any of the equations
\begin{equation}\label{eqs}
\d_t u - \Delta u = |u|^p, \quad\quad\quad \d_t^2 u + a\d_t u -\Delta u = |u|^p,
\end{equation}
where $a\in \R$.  More precisely, the conclusions of Theorems \ref{wave} and \ref{heat} (both (a) and (b)) remain valid.
That is, the instability occurs by blow up in finite time.

If $p\geq p_c$, on the other hand, then $\phi$ is a linearly stable solution in the sense that the linearized operator
has no negative spectrum.
\end{theorem}

\begin{proof}
Let $Q$ denote the quadratic in \eqref{Q1} and consider the expression
\begin{equation}\label{Qn1}
\mathcal{Q}(p) \equiv 4(p-1)^2 \cdot \left[ \left( \frac{n-2}{2} \right)^2- p\cdot Q \left( \frac{2}{p-1} \right)\right].
\end{equation}
When this is non-negative, Lemma \ref{sp} implies that $\phi$ is linearly stable.  When it is negative, on the other
hand, Lemma \ref{ip} implies that the linearized operator has a negative eigenvalue. Since $\phi$ vanishes at infinity,
each of the functions
\begin{equation*}
\phi,\quad\quad f(\phi) = \phi^p, \quad\quad f'(\phi)=p\phi^{p-1}
\end{equation*}
is bounded, whence (A3)-(A6) hold and the linearized operator has essential spectrum
\begin{equation*}
\sigma_{\text{ess}}(-\Delta-p\phi^{p-1}) = \sigma_{\text{ess}}(-\Delta) = [0,\infty),
\end{equation*}
whence (A1)-(A2) follow by standard arguments, as in the proof of Theorem \ref{new}.

In view of these observations, it thus suffices to check that
\begin{equation}\label{goal}
\mathcal{Q}(p) <0 \quad\iff\quad \frac{n+2}{n-2} \leq p< p_c.
\end{equation}
Combining our definitions \eqref{Q1} and \eqref{Qn1}, we have
\begin{align}\label{Qp}
\mathcal{Q}(p) = (n-2)(n-10) \cdot p^2 - 2(n^2-8n+4)\cdot p + (n-2)^2.
\end{align}
The last equation easily leads to
\begin{equation}\label{Qpn}
\mathcal{Q}\left( \frac{n+2}{n-2} \right) = -\frac{64}{n-2} < 0.
\end{equation}

\Case{1}
Suppose first that $n>10$.  Then $\mathcal{Q}(p)$ has two roots, the largest one of which is $p_c$, given by \eqref{pc}.
Since $\mathcal{Q}(p)$ is negative only between its two roots, \eqref{Qpn} implies that $\mathcal{Q}(p)$ is negative if
$\frac{n+2}{n-2} \leq p< p_c$ and non-negative if $p\geq p_c$.

\Case{2} Suppose now that $2<n\leq 10$.  Then the computation
\begin{align*}
\mathcal{Q}'(p) &= 2(n-2)(n-10) \left(p - \frac{n+2}{n-2} \right) - 48
\end{align*}
forces $\mathcal{Q}(p)$ to be decreasing for each $p\geq \frac{n+2}{n-2}$, and hence negative by \eqref{Qpn}.  This shows
that the desired condition \eqref{goal} holds with $p_c=\infty$.
\end{proof}

\begin{rem}\label{re-h}
When it comes to the nonlinear heat equation \eqref{eqs}, the above result of Theorem \ref{old} is originally due to Gui,
Ni and Wang \cite{GNW}.  Their proof is based on an earlier result of Wang \cite{Wang}, according to which the graphs of
two steady states do not intersect if $p\geq p_c$, while the graphs of any two steady states do intersect if $p<p_c$.
Using this characterization of the critical value $p_c$, they are then able to employ comparison arguments which rely on
the strong maximum principle. As far as the unstable case $p<p_c$ is concerned, they show that blow up occurs provided
the perturbation $\psi_0$ in \eqref{3} satisfies $0\leq \psi_0\nequiv 0$, while solutions exist globally provided  $0\geq
\psi_0\nequiv 0$.  Theorem \ref{old} yields a refinement of their blow up result for the nonlinear heat equation since
our exact assumption \eqref{per2} on $\psi_0$ reads
\begin{equation*}
\int_{\R^n} \chi(x) \,\psi_0(x) \:dx > 0
\end{equation*}
for a certain positive function $\chi(x)$. As for the stable case $p\geq p_c$, the result in \cite{GNW} is much stronger
than Theorem \ref{old} because it proves nonlinear stability, not merely linear stability.  However, a proof of nonlinear
stability for the analogous hyperbolic problem is difficult to imagine.
\end{rem}

\begin{rem}\label{re-w}
When it comes to the wave equation \eqref{eqs} in the critical case $p=\frac{n+2}{n-2}$, our conclusions are closely
related to a result of Krieger and Schlag \cite{KS} in three dimensions. According to Theorem \ref{ex}, the equation
\begin{equation}\label{ex5}
\d_t^2 u - \Delta u = |u|^5, \quad\quad x\in\R^3
\end{equation}
admits a family of steady states
\begin{equation}\label{curve}
\phi_\la(x) = \frac{(3\la^2)^{1/4}}{\sqrt{\la^2 + |x|^2}} \:,\quad\quad \la>0,
\end{equation}
all of which are unstable by Theorem \ref{old}. For the equation with $u^5$ rather than $|u|^5$, Krieger and Schlag
\cite{KS} construct a stable manifold $\Sigma$ associated with $\phi_1$.  That is, if the initial data
\begin{equation*}
u(x,0)= \phi_1(x)+\psi_0(x), \quad\quad \d_t u(x,0)= \psi_1(x)
\end{equation*}
correspond to a small perturbation $(\psi_0,\psi_1)\in \Sigma$, then the solution to \eqref{ex5} exists globally and
remains near the curve $\{\phi_\la\}_{\la>0}$ for all times.  Numerical computations \cite{Num1, Num2} suggest that the
stable manifold $\Sigma$ ought to represent the borderline case between global existence and blow up of solutions.  To
see that Theorem \ref{old} provides a partial proof of this conjecture, we note that the tangent plane to $\Sigma$ at the
origin is given by the equation
\begin{equation}\label{per3}
\sigma \int_{\R^n} \chi(x)\,\psi_0(x) \:dx + \int_{\R^n} \chi(x)\,\psi_1(x) \:dx =0,
\end{equation}
where $\sigma>0$ and $-\sigma^2$ is the first eigenvalue of the linearized operator; see \cite{KS}.  Thus the set of
perturbations \eqref{per1} to which our instability theorem applies is precisely the set of perturbations that lie
strictly above the tangent plane to $\Sigma$ at the origin.
\end{rem}

\end{document}